\documentclass{amsart}

\usepackage{amssymb}
\usepackage{amsmath}
\usepackage{amsthm}
\usepackage{amsfonts}
\usepackage{amstext}
\usepackage{algorithmic}
\usepackage{algorithm}
\usepackage{graphicx}
\usepackage{epstopdf}
\usepackage[all]{xy}

\allowdisplaybreaks

\newtheorem{theorem}{Theorem}[section]
\newtheorem{lemma}[theorem]{Lemma}
\newtheorem{proposition}[theorem]{Proposition}
\newtheorem{corollary}[theorem]{Corollary}
\theoremstyle{definition}
\newtheorem{definition}[theorem]{Definition}
\newtheorem{example}[theorem]{Example}
\theoremstyle{remark}
\newtheorem{remark}[theorem]{\bf{Remark}}

\renewcommand{\O}{{\mathbb{O}}}

\newcommand{\Z}{{\mathbb{Z}}}

\newcommand{\Rcal}{{\mathcal{R}}}

\newcommand{\tens}{\otimes}
\newcommand{\id}{{\rm id}}

\renewcommand{\o}{{}_{(1)}}
\renewcommand{\t}{{}_{(2)}}

\newcommand{\eps}{\epsilon}
\newcommand{\vareps}{\varepsilon}
\newcommand{\la}{{\triangleright}}
\newcommand{\ra}{{\triangleleft}}

\newcommand{\bos}{{> \!\!\!\triangleleft\kern-.33em\cdot\, }}
\newcommand{\crosscop}{\,_{\Rcal}\!\!\blacktriangleright\kern -4 pt\blacktriangleleft}

\newcommand{\lcocross}{{>\!\!\blacktriangleleft}}
\newcommand{\bicross}{\triangleright\!\!\!\blacktriangleleft}
\newcommand{\lrbicross}{\blacktriangleright\!\!\!\triangleleft}

\begin{document}

\title{Bicrossproduct Hopf quasigroups}
\date{27th October 2009}
\author{Jennifer Klim \& Shahn Majid}
\address{School of Mathematical Sciences\\ Queen Mary University of London\\ 327 Mile End Rd, London E1 4NS, UK}

\begin{abstract}
We recall the notion of Hopf quasigroups introduced in \cite{Klim09}. We construct  a bicrossproduct Hopf quasigroup $kM\bicross k(G)$ from every group $X$ with a finite subgroup $G\subset X$ and IP quasigroup transversal $M\subset X$  subject to certain conditions. We identify the octonions quasigroup $G_\O$ as transversal in an order 128 group $X$ with subgroup $\Z_2^3$ and hence obtain a Hopf quasigroup $kG_{\O}\lcocross k(\Z_2^3)$ as a particular case of our construction.
\end{abstract}

\maketitle

\section{Introduction}

In \cite{Klim09} we defined the notion of a Hopf quasigroup as a non-associative algebra with a coassociative coproduct together with an antipode $S$ satisfying certain linearised IP quasigroup identities. In this paper we find a further class of noncommutative and noncocommutative examples by means of a bicrossproduct construction. Such a construction was successfully used to construct the first Hopf algebras associated to group factorisations and we do the same now for a factorisation of a group $X$ into a finite subgroup $G$ and a transversal $M$ containing the identity. Among the conditions we find for the antipode is that the product on $M$ defined by the product of coset representatives modulo $G$ makes $M$ a IP quasigroup.  We give necessary and sufficient conditions within our construction. Throughout the paper unless qualified otherwise a quasigroup will  mean with two-sided inverse property (also called an IP loop).

An outline of the paper is as follows; in Section 2 we recall the definition of a Hopf quasigroup and briefly mention without proof some of the results in \cite{Klim09}.   In Section 3, we begin with the matched pair approach to subgroups and transversals \cite{Beggs03}. It is a standard result in quasigroup theory that if $X$ is a group, $G\subset X$ a subgroup and $M$ a set of left coset representatives (a transversal) then $M$ acquires a left quasigroup structure, and variously more structure as $M$ is variously restricted\cite{Smi:book}. In the matched pair approach this data is developed in terms of  a right action $\ra$ of $G$ on $M$ and a map $\la: M\times G\to G$ which is some kind of left `quasi' action, a `cocycle' $\tau:M\times M\to G$  and an induced (generally nonassociative) product $\cdot$ on $M$. We extend \cite{Beggs03} to an analysis of when $M$ is an IP quasigroup. 

We are then in position in Section 4 to apply a process of semidualization in which one factor is dualised to obtain a  bicrossproduct $kM\bicross k(G)$ as a semidirect product by $\la$ and semidirect coproduct by $\ra$. We show that this is a Hopf quasigroup iff $M$ is a quasigroup and certain identities hold for $\tau$.  Dualising the other factor gives a Hopf coquasigroup $k(M)\lrbicross kG$. We require $G$ respectively $M$ to be finite.
We also provide an example where $X=\Z_2^3\ltimes Cl_3$ an order 128 non-Abelian group, $G=\Z_2^3$, and $M=G_{\O}$, the octonion quasigroup. Here $\la$ is trivial, and we obtain a Hopf quasigroup $kG_{\O}\lcocross k(G)$. Unless stated otherwise, we work over a general field $k$.

\section{Hopf Quasigroups}

We recall that an (inverse property) {\em quasigroup} is a set $M$ with a product, denoted for the moment by omission, an identity $e$ and for each $s\in M$ an element $s^{-1}\in M$ such that 
\[ s^{-1}(s t)=t,\quad (t s)s^{-1}=t,\quad \forall t\in M.\]
A quasigroup  is {\em Moufang} if $s(t(sr))=((st)s)r$ for all $s,t,r\in M$.

In \cite{Klim09} we linearised these notions to Hopf quasigroup and Moufang Hopf quasigroup in the same way that a Hopf algebra linearises the notion of a group.

\begin{definition}\cite{Klim09}  A \textit{Hopf quasigroup}  is a possibly-nonassociative but unital algebra $H$ equipped with algebra homomorphisms $\Delta:H\to H\tens H$, $\eps:H\to k$ forming a coassociative coalgebra and a map $S:H\to H$ such that
\begin{equation*} m(\id\tens m)(S\tens\id\tens\id)(\Delta\tens\id)=\eps\tens \id = m(\id\tens m)(\id\tens S\tens\id)(\Delta\tens\id) \end{equation*}
\begin{equation*} m(m\tens\id)(\id\tens S\tens\id)(\id\tens\Delta)=\id\tens\eps = m(m\tens\id)(\id\tens\id\tens S)(\id\tens\Delta) \end{equation*}
\end{definition}
One can write these more explicitly as
\[ \sum Sh\o (h\t g)=\sum h\o(Sh\t g)=\sum (g Sh\o)h\t=\sum (g h\o)Sh\t=\eps(h)g\]
for all $h,g\in H$, where we write $\Delta h=\sum h\o\tens h\t$. In this notation the Hopf quasigroup $H$ is called \textit{Moufang} if 
\[\sum h\o(g(h\t f))=\sum ((h\o g)h\t)f \quad \forall h,g,f\in H\]

The conditions in Definition~2.1 are stronger than the usual Hopf algebra antipode axioms and compensate  for $H$ nonassociative.  The paper \cite{Klim09} provides the first results establishing a  viable Hopf-like theory. For example, $S$ is antimultiplicative and anticomultiplicative in the sense
\[ S(hg)=(Sg)(Sh),\quad \Delta(S(h))=S(h\t)\tens S(h\o)\]
for all $h,g\in H$.  Clearly an actual (inverse property) quasigroup $M$ linearizes to a Hopf  quasigroup algebra $kM$ with grouplike coproduct on elements of $M$ and linear extension of the product and inverse, and is Moufang if $M$ is. As with Hopf algebras, the theory unifies group and enveloping algebras:

\begin{proposition}\cite{Klim09} For $(L,[\ ,\ ])$ a Malt'sev algebra over $k$ not of characteristic 2,3,  the enveloping algebra $U(L)$  in \cite{SS} is a Moufang Hopf quasigroup.
\end{proposition}

Also as in Hopf algebra theory, we have a dual notion, which we called a \textit{Hopf coquasigroup} by reversing the arrows on each of the maps. 

\begin{proposition}\cite{Klim09}
The algebraic variety $k[S^7]$ is a Moufang Hopf coquasigroup.
\end{proposition}

We used $S^7\subset \O$ and the cochain twist approach to the octonions of \cite{AlbMa}. In this context there is a natural action of $\Z_2^3$ on $k[S^7]$, which leads\cite{Klim09} to a cross product $k[S^7]\rtimes \Z_2^n$ as a first example of a noncommutative noncocommutative Hopf coquasigroup. The dual $k(G_\O)\rtimes \Z_2^3$ of the bicrossproduct Hopf quasigroup now obtained in Section~4  could be viewed as a quotient of this corresponding to the inclusion $G_\O\subset S^7$.

\section{Matched pair approach to coset quasigroups}

Let $X$ be a group and $G\subset X$ be a subgroup. A transversal $M\subset X$ means a set of left coset representatives  of $G$, i.e. such that for every $x\in X$ there exists a unique $s\in M$ such that $x\in Gs$. We follow the approach to this data in \cite{Beggs03} since this is closest to the more familiar `matched pair of actions' in the group factorisation case leading to ordinary Hopf algebra  bicrossproducts\cite{Ma:book}. We assume throughout that $e\in M$, where $e$ is the identity of $X$. Let $s,t\in M$, $u\in G$ and define a cocycle $\tau:M\times M\to G$ and a product $\cdot$ on $M$ by
\begin{equation}\label{crossrel} st = \tau(s,t)(s\cdot t),\quad  su=(s\la u)(s\ra u) \end{equation}
for some kind of actions $\la$ of $M$ on $G$, and $\ra$ of $G$ on $M$. Most of the following identities  are obtained from the associativity of $X$ on various combinations of elements using the unique factorisation.

\begin{proposition}\cite{Beggs03} The following identities between $(M,\cdot)$ and $\la,\ra,\tau$ hold for all $s,t,r\in M$ and $u,v\in G$:
\begin{equation} \tau(s,t)\tau(s\cdot t,r) = (s\la\tau(t,r))\tau(s\ra\tau(t,r),t\cdot r) \label{2} \end{equation}
\begin{equation} s\ra(uv) = (s\ra u)\ra v \label{3}\end{equation}
\begin{equation} (s\cdot t)\ra u = (s\ra(t\la u))\cdot (t\ra u) \label{4}\end{equation}
\begin{equation} (s\cdot t)\cdot r = (s\ra\tau(t,r))\cdot(t\cdot r)  \label{5} \end{equation}
\begin{equation} \tau(s,t)((s\cdot t)\la u) = (s\la(t\la u))\tau(s\ra(t\la u),t\ra u) \label{6} \end{equation}
\begin{equation} s\la(uv) = (s\la u)((s\ra u)\la v) \label{7} \end{equation}
\begin{equation} s\la e = e,\quad s\ra e = s,\quad \tau(s,e)=e=\tau(e,s),\quad s\cdot e= s = e\cdot s.\label{8} \end{equation}
\end{proposition}

Let $s\in M$, if $s^{-1}$ in $X$ lies in $M$ then
\[ e= s^{-1}s=\tau(s^{-1},s)s^{-1}\cdot s  \]
\[ e = ss^{-1}=\tau(s,s^{-1})s\cdot s^{-1} \]
and we find $\tau(s^{-1},s)=e =\tau(s,s^{-1})$ and $s^{-1}\cdot s = e = s\cdot s^{-1}$.  More generally, we have:

\begin{lemma}\label{lemma3.2}
Let $s\in M$. There is a unique $s^{-L}$ such that $s^{-1}=\tau^{-1}(s^{-L},s)s^{-L}$. We have
\begin{equation}s^{-L}\cdot s	= e\label{SL1}\end{equation}
\begin{equation} (s\ra\tau^{-1}(s^{-L},s))\cdot s^{-L} = e \label{SL2}\end{equation}
\begin{equation} s\la \tau^{-1}(s^{-L},s)	= \tau^{-1}(s\ra\tau^{-1}(s^{-L},s),s^{-L})\label{SL3}\end{equation}
\end{lemma}
\proof We let $s^{-1}=\chi(s)s^-$ where $\chi(s)\in G,s^-\in M$. Then 
\[ e = s^{-1}s = \chi(s)s^-s = \chi(s)\tau^{-1}(s^-,s)s^-\cdot s \]
So we find $\chi(s)=\tau^{-1}(s^-,s)$ and $s^-\cdot s=e$. Since $s^-$ is a left inverse of $s$ in $M$, we shall label it $s^{-L}$. Now, using (\ref{crossrel}) we have
\begin{eqnarray*}
e &	=	& ss^{-1} = s\tau^{-1}(s^{-L},s)s^{-L} = (s\la\tau^{-1}(s^{-L},s))(s\ra\tau^{-1}(s^{-L},s))s^{-L} \\
	&= & (s\la\tau^{-1}(s^{-L},s))\tau(s\ra\tau^{-1}(s^{-L},s),s^{-L})(s\ra\tau^{-1}(s^{-L},s))\cdot s^{-L}
\end{eqnarray*}
which implies $s\la\tau^{-1}(s^{-L},s) = \tau^{-1}(s\ra\tau^{-1}(s^{-L},s),s^{-L})$ and $(s\ra\tau^{-1}(s^{-L},s))\cdot s^{-L} = e$. \endproof

It is a familiar fact\cite{Smi:book} that  the transversal $M$ under our assumptions has the structure of a right quasigroup (with identity in our case, i.e. a right loop). This means the existence of a `division' map $/:M\times M\to M$ such that $(t\cdot s)/s=t=(t/s)\cdot s$ for all $s,t\in M$ and means in particular that we have right cancellation: $t\cdot s=t'\cdot s\Rightarrow t=t'$. This is the content of Lemma~\ref{lemma3.2} in terms of the matched pair data $\la,\ra,\tau$:

\begin{proposition} $M$ in the setting of  Proposition~3.1 is a right quasigroup (with identity). The division map is
\[ t/s=t\ra\tau^{-1}(s^{-L},s)\cdot s^{-L}.\]
Moreover,
\[ (s^{-L})^{-L}=s\ra\tau^{-1}(s^{-L},s).\]
\end{proposition}
\proof $(t/s)\cdot s=((t\ra\tau^{-1}(s^{-L},s))\cdot s^{-L})\cdot s=t\cdot(s^{-L}\cdot s)=t$ by (\ref{5}) and (\ref{SL1}). On the other side 
\begin{eqnarray*}
(t\cdot s)/s &=& ((t\cdot s)\ra\tau^{-1}(s^{-L},s))\cdot s^{-L}=((t\ra(s\la\tau^{-1}(s^{-L},s)))\cdot (s\ra\tau^{-1}(s^{-L},s)))\cdot s^{-L}\\
&=&((t\ra\tau^{-1}(s\ra \tau^{-1}(s^{-L},s),s^{-L}))\cdot (s\ra\tau^{-1}(s^{-L},s)))\cdot s^{-L}\\
&=&t\cdot ( (s\ra\tau^{-1}(s^{-L},s))\cdot s^{-L} )=t\end{eqnarray*}
using (\ref{4}), (\ref{SL3}), (\ref{5}) and (\ref{SL2}). Once we have right cancellation, we know that $s^{-L}$ is the unique left inverse for each $s\in M$, hence (\ref{SL2}) implies the second part of the proposition.
\endproof

Using (\ref{4}) and (\ref{7}) respectively, we also obtain the following useful identities:
\begin{equation} (t\la v)^{-1} = (t\ra v)\la v^{-1} \label{12}\end{equation}
\begin{equation} (t\ra v)^{-L} = t^{-L}\ra(t\la v)\label{13} \end{equation}

Finally, we will know that we have captured {\em all} of the input data into Proposition~3.1 if we can rebuild $X$ from $\la,\ra$ and $\tau$ and the identities there. In the group factorisation case this is the construction of a double cross product group $G\bowtie M$ from a matched pair data. Most of this is in  \cite{Beggs03}. 

\begin{proposition}
Suppose that $G$ is a group and the maps $\la,\ra$ and $\tau$ satisfy the identities (\ref{2})-(\ref{7}). Then the set  $G\times M$ acquires the structure of a group $G\bowtie M$ by 
\[ (u,s)(v,t) = (u(s\la v)\tau(s\ra v,t),(s\ra v)\cdot t), \quad (u,s)^{-1} = (\tau^{-1}(s^{-L},s)(s^{-L}\la u^{-1}),s^{-L}\ra u^{-1})\]
and identity $(e,e)$. This group is isomorphic to $X$ with subgroup $(G,e)$ and transversal $(e,M)$.
\label{Beggs}
\end{proposition}
\proof
The required form of the product is easily obtained from (\ref{crossrel}) and $(u,s)(v,t)=(u,e)(e,s)(v,e)(e,t)$ after which everything can be verified by direct computation. We only need to verify that we have a right inverse as the rest was covered in \cite{Beggs03}. We compute
\begin{eqnarray*}
&&\kern-20pt (u,s)(\tau^{-1}(s^{-L},s)(s^{-L}\la u^{-1}),s^{-L}\ra u^{-1})\\
	&	=	& (u(s\la\tau^{-1}(s^{-L},s)(s^{-L}\la u^{-1}))\tau(s\ra(\tau^{-1}(s^{-L},s)(s^{-L}\la u^{-1})),s^{-L}\ra u^{-1}),\\
	&		& \quad \quad (s\ra(\tau^{-1}(s^{-L},s)(s^{-L}\la u^{-1})))\cdot(s^{-L}\ra u^{-1}))\\
	&	=	& (u(s\la\tau^{-1}(s^{-L},s))((s\ra\tau^{-1}(s^{-L},s))\la(s^{-L}\la u^{-1}))\tau((s\ra\tau^{-1}(s^{-L},s))\ra(s^{-L}\la u^{-1}),s^{-L}\ra u^{-1}),\\
	&		&\quad\quad ((s\ra\tau^{-1}(s^{-L},s))\ra(s^{-L}\la u^{-1}))\cdot (s^{-L}\ra u^{-1}))\\
	&	=	& (u(s\la\tau^{-1}(s^{-L},s))\tau(s\ra\tau^{-1}(s^{-L},s),s^{-L})(((s\ra\tau^{-1}(s^{-L},s))\cdot s^{-L})\la u^{-1}),\\
	&		&\quad\quad ((s\ra\tau^{-1}(s^{-L},s))\cdot s^{-L})\ra u^{-1})\\
	&	=	& (u(s\la\tau^{-1}(s^{-L},s))\tau(s\ra\tau^{-1}(s^{-L},s),s^{-L})(e\la u^{-1}),e\ra u^{-1})\\
	&	=	& (u(s\la\tau^{-1}(s^{-L},s))\tau(s\ra\tau^{-1}(s^{-L},s),s^{-L})u^{-1},e)=(e,e)
\end{eqnarray*}
The second equality uses (\ref{3}) and (\ref{7}), the third uses (\ref{4}) and (\ref{6}), the fourth uses (\ref{SL2}), and the final equality uses (\ref{SL3})  in Lemma~\ref{lemma3.2}. \endproof

Continuing our analysis, depending on $M$ we may also have a right inverse $s^{-R}$ of $s$ in $M$. If so we have the following compatibility relations 
\[ s^{-L}\la\tau(s,s^{-R}) = \tau(s^{-L},s) \]
\[ s^{-L} \ra \tau(s,s^{-R})=s^{-R} \]
by considering $(s^{-L}\cdot s)\cdot s^{-R} = s^{-L}\ra\tau(s,s^{-R})\cdot (s\cdot s^{-R})$ and similar. 
The two inverses will not generally coincide unless $M$ is a quasigroup.

\begin{proposition}
In the setting of Proposition~3.1, the following  are equivalent
\begin{enumerate}
	\item[(i)] $\forall s\in M$, $s\la G=G$
 \item[(ii)] $\forall s\in M$ there exists a right inverse $s^{-R}\in M$
	\item[(iii)] $X=MG.$
\end{enumerate}\label{prop3.5}
\end{proposition}
\proof We will show that $(i)\Rightarrow(iii)\Rightarrow(ii)\Rightarrow(i)$.

Suppose $u\in G$, $s\in M$ and that (i) holds. We seek $t,v$ such that $us=tv=(t\la v)(t\ra v)$. Hence if we can find $v$ we will have a unique $t=s\ra v^{-1}$. It remains to find $v$ solving $(s\ra v^{-1})\la v=u$. By (\ref{12}) this is equivalent to $(s\la v^{-1})^{-1}=u$ or $s\la v^{-1}=u^{-1}$. Such $v$ exists under our assumption,  hence (iii) holds. We see also that if $s\la(\ )$ is bijective then the factorisation as $MG$ is unique.

Suppose (iii) holds, then for $s\in M$,  $s^{-1}=s^-\psi(s)$ for some $s^-\in M,\psi(s)\in G$.
\[ e=ss^{-1}=ss^-\psi(s) = \tau(s,s^-)(s\cdot s^-)\psi(s) = \tau(s,s^-)((s\cdot s^-)\la \psi(s))((s\cdot s^-)\ra \psi(s))\]
so $(s\cdot s^-)\ra \psi(s)=e\Rightarrow s\cdot s^-=e$, implying right inverses $s^{-R}$ exist, and (ii) holds. We also learn that $\psi(s)=\tau(s,s^{-R})$.  Note that if the factorisation as $MG$ is unique we can write $st=s\circ t\sigma(s,t)\in MG$ as a left-right reversal of our previous theory. Then 
\[ss^{-R}=\tau(s,s^{-R})s\cdot s^{-R}=s\circ s^{-R}\sigma(s,s^{-R})=((s\circ s^{-R})\la\sigma(s,s^{-R}))((s\circ s^{-R})\ra\sigma(s,s^{-R}))\]
from which we see that $e=(s\circ s^{-R})\ra\sigma(s,s^{-R})$ and hence $s\circ s^{-R}=e$. But as $(M,\circ)$ has left cancellation, $s^{-R}$ is uniquely determined. 

Finally, suppose (ii) holds and let $s\in M,u\in G$. If there exists $v\in G$ such that $s\la v=u$ then
\begin{eqnarray*}
v	&	=	& (s^{-L}\cdot s)\la v\\
	&	=	& \tau^{-1}(s^{-L},s)s^{-L}\la(s\la v)\tau(s^{-L}\ra(s\la v),s\ra v)\\
	&	=	& \tau^{-1}(s^{-L},s)s^{-L}\la u\tau(s^{-L}\ra u,s\ra v)
\end{eqnarray*}
by (\ref{6}) and $(s\ra v)^{-L}=s^{-L}\ra u$ by (\ref{13}), which can be solved by $s\ra v = (s^{-L}\ra u)^{-R}$. Accordingly, we take
\[ v=\tau^{-1}(s^{-L},s)(s^{-L}\la u)\tau(s^{-L}\ra u,(s^{-L}\ra u)^{-R}) \]
 as a definition and verify,
\begin{eqnarray*}
\lefteqn{s\la(\tau^{-1}(s^{-L},s)(s^{-L}\la u)\tau(s^{-L}\ra u,(s^{-L}\ra u)^{-R}))}\\
	&	=	& (s\la\tau^{-1}(s^{-L},s))((s\ra\tau^{-1}(s^{-L},s))\la((s^{-L}\la u)\tau(s^{-L}\ra u,(s^{-L}\ra u)^{-R}))\\
	&	=	& (s\la\tau^{-1}(s^{-L},s))((s\ra\tau^{-1}(s^{-L},s))\la(s^{-L}\la u))\\
	&		& \quad\quad (((s\ra\tau^{-1}(s^{-L},s))\ra(s^{-L}\la u))\la\tau(s^{-L}\ra u,(s^{-L}\ra u)^{-R}))\\
	&	=	& (s\la\tau^{-1}(s^{-L},s))\tau(s\ra\tau^{-1}(s^{-L},s),s^{-L})((s\ra\tau^{-1}(s^{-L},s))\cdot s^{-L}\la u)\\
	&		& \quad\quad \tau^{-1}((s\ra\tau^{-1}(s^{-L},s))\ra(s^{-L}\la u),s^{-L}\ra u)\\
	&		& \quad\quad\quad(((s\ra\tau^{-1}(s^{-L},s))\ra(s^{-L}\la u))\la\tau(s^{-L}\ra u,(s^{-L}\ra u)^{-R}))\\
	&	=	& (e\la u)\tau^{-1}((s\ra\tau^{-1}(s^{-L},s))\ra(s^{-L}\la u),s^{-L}\ra u)\\
	&		& \quad\quad\quad(((s\ra\tau^{-1}(s^{-L},s))\ra(s^{-L}\la u))\la\tau(s^{-L}\ra u,(s^{-L}\ra u)^{-R}))\\
	&	=	& u\tau(((s\ra\tau^{-1}(s^{-L},s))\ra(s^{-L}\la u))\cdot(s^{-L}\ra u),(s^{-L}\ra u)^{-R})\\
	&		& \quad\quad \tau^{-1}(((s\ra\tau^{-1}(s^{-L},s))\ra(s^{-L}\la u))\ra\tau(s^{-L}\ra u,(s^{-L}\ra u)^{-R}),(s^{-L}\ra u)\cdot(s^{-L}\ra u)^{-R})\\
	&	=	& u\tau((s\ra\tau^{-1}(s^{-L},s))\cdot s^{-L}\ra u,(s^{-L}\ra u)^{-R})\\
	&		& \quad\quad \tau^{-1}(((s\ra\tau^{-1}(s^{-L},s))\ra(s^{-L}\la u))\ra\tau(s^{-L}\ra u,(s^{-L}\ra u)^{-R}),e)\\
	&	=	& u\tau(e\ra u,(s^{-L}\ra u)^{-R})\\
	&	=	& u
\end{eqnarray*}
as required. Here the first and second equalities use (\ref{7}), the third uses (\ref{6}). The fourth equality uses (\ref{SL2}) and (\ref{SL3}), the fifth uses (\ref{2}), the sixth uses (\ref{4}), and the seventh uses (\ref{SL2}). We also see that if the right inverse is unique then so is $v$ because $s\ra v$ and hence $v$ then necessarily have the form used. \endproof

Note that if $G$ is finite then $s\la(\ )$ in condition (i) will  be bijective. We have noted in the proof that this is equivalent to the factorisation in (iii) being unique and to the the right inverse in (ii) being unique.

\begin{proposition}\label{prop3.6}
$M$ in the setting of Proposition~3.1  is a (two-sided inverse property) quasigroup {\em iff} 
\[ t = t\ra\tau(s^{-L},s) \]
\[  s^{-L} = s^{-L}\ra\tau(s,t)\]
for all $s,t\in M$. In this case $(s^{-L})^{-L} = s$, i.e. $s^{-L}=s^{-R}$.
\end{proposition}
\proof
Suppose $M$ is an inverse property quasigroup, then for all $s,t\in M$,
\[ t = (t\cdot s^{-L})\cdot s = (t\ra\tau(s^{-L},s))\cdot(s^{-L}\cdot s)=(t\ra\tau(s^{-L},s))\cdot e= t\ra\tau(s^{-L},s). \]
Similarly, using (\ref{5}) we have $t=s^{-L}\cdot(s\cdot t)=((s^{-L}\ra\tau^{-1}(s,t))\cdot s)\cdot t$ and conclude that $s^{-L}\ra\tau(s,t) = s^{-L}$  (and setting $t=s^{-R}$  in this confirms that $s^{-L}=s^{-R}$ as it must in an IP quasigroup.)

 Conversely, suppose $t=t\ra \tau(s^{-L},s)$ then from Proposition~3.3 we see that $t/s=t\cdot s^{-L}$ and hence $(t\cdot s)\cdot s^{-L}=t$.  For the other side if we suppose that $s^{-L} = s^{-L}\ra\tau(s,t)$ then $s^{-L}\cdot(s\cdot t)=((s^{-L}\ra\tau^{-1}(s,t)) \cdot s)\cdot t =(s^{-L}\cdot s)\cdot t=t$ using (\ref{5}). \endproof

When $M$ is a quasigroup, we shall continue to denote the (left and right) inverse of $s\in M$ by $s^{-L}$ to distinguish it from $s^{-1}\in X$. We shall need two further elementary lemmas.
 
\begin{lemma} In the setting of Proposition~3.1, 
\[  (s\ra u)^{-L}\la(s\la u)^{-1} = u^{-1} \]
holds  for  $s\in M$, $u\in G$ {\em iff}   $u^{-1}\tau(s^{-L},s)u=\tau(s^{-L}\ra(s\la u),s\ra u)$.
\label{lemma3.7}
\end{lemma}
\proof
\begin{eqnarray*}
(s\ra u)^{-L}\la(s\la u)^{-1}	&	=	& (s^{-L}\ra(s\la u))\la(s\la u)^{-1} \\
	&	=	&	(s^{-L}\la(s\la u))^{-1}\\
	&	=	&	\tau(s^{-L}\ra(s\la u),s\ra u)u^{-1}\tau^{-1}(s^{-L},s)\\
	\end{eqnarray*}
using (\ref{13}),  (\ref{12}) and (\ref{6}). So equality of this to $u^{-1}$ is precisely the stated condition for $\tau$.  From the proof we see that the latter is also equivalent to $s^{-L}\la(s\la u)=u$. \endproof

\begin{lemma}\label{lemma3.8} In the setting of Proposition~3.1, we have
\[ (s\ra u)^{-L}\ra(s\la u)^{-1} = s^{-L} \]
\label{lemma2} for  all $s\in M$ and $u\in G$. \end{lemma}
\proof Using (\ref{12}) and (\ref{13}),
\[ (s\ra u)^{-L}\ra(s\la u)^{-1} = (s^{-L}\ra(s\la u))\ra(s\la u)^{-1} = s^{-L}\ra ((s\la u)(s\la u)^{-1}) = s^{-L} \]
\endproof

In summary, we have dissected the usual coset construction in terms of properties of the matched pair data $\la,\ra,\tau$. We will now use this data to construct something different.

\section{Bicrossproduct $kM \bicross k(G)$}

With $M$ and $G$ as in Proposition~3.1 and $G$ finite, we consider the bicrossproduct $kM\bicross k(G)$. By $kM$ we mean the vector space on $M$ with its product $\cdot$ extended linearly and a group-like coproduct on basis elements. By $k(G)$ we mean the functions on $G$ with  usual pointwise multiplication and standard coproduct.  From the `quasi-action' $\la$ in our matched pair data we have a right action on $k(G)$ and make a cross product algebra. From the action $\ra$ in our data we have a left coaction of $k(G)$ on $kM$ and make cross coproduct coalgebra. This gives  $kM\bicross k(G)$ with exactly the same formulae as for Hopf algebra bicrossproducts. Following the conventions of \cite{Ma:book} we take basis $\{s\tens\delta_u|s\in M,u\in G\}$ and the algebra and coalgebra structure are explicitly
\[ (s\tens\delta_u)(t\tens\delta_v) = \delta_{u,t\la v} (s\cdot t\tens\delta_v) \]
\[ 1 =\sum_u e\tens\delta_u \]
\[ \Delta(s\tens\delta_u) = \sum_{ab=u} (s\tens\delta_a)\tens(s\ra a\tens\delta_b) \]
\[ \vareps(s\tens\delta_u) = \delta_{u,e} \]

Notice that since $\ra$ is an actual group action it induces an actual coaction of $k(G)$ and $\Delta$ is therefore coassociative as a standard cross coproduct on the coalgebra of $kM$. Put it another way, the proof is identical to the proof \cite{Ma:book} for ordinary bicrossproducts; it is not sensitive to $M$ not necessarily being a group. Likewise when checking the homomorphism property  $\Delta((s\tens\delta_u)\cdot (t\tens\delta_v))=(\Delta(s\tens\delta_u))\cdot(\Delta(t\tens\delta_v))$ we only ever encounter in each tensor factor one product in $M$ and one application of $\ra$, so we never use any of the identities in Proposition~3.1 involving $\tau$. All the others are identical to those in the matched pair conditions for ordinary bicrossproducts, so the proof  is again line by line identical. Incidentally, this proof has a nice diagrammatic formulation in terms of subdividing squares\cite{Ma:book}. Hence the only issue is the antipode.

\begin{lemma}\label{lemma4.1} Suppose $M$ in Proposition~3.1 has right inverses and $G$ is finite. The antipode of the bicrossproduct $kM\bicross k(G)$, if it is a Hopf quasigroup, necessarily takes the form
\[ S(s\tens\delta_u) = (s\ra u)^{-L}\tens\delta_{(s\la u)^{-1}} \]
for all $s\in M$ and $u\in G$. \end{lemma}
\proof
Consider $S(s\tens\delta_u) = \sum_t t\tens f_t^{u,s}$ for some functions with  $f_t^{u,s}\in k(G)$.  Then considering the usual antipode property alone (set $g=1$ in the explicit formulae after Definition~2.1) applied to $h=s\tens \delta_{v}$ we have 
\begin{eqnarray*}
e\tens 1\delta_{v,e}&	=	& \sum (S(s\tens \delta_v)\o)(s\tens \delta_v)\t	\\
	&	=	& \sum_{a} S(s\tens\delta_a)(s\ra a\tens\delta_{a^{-1}v})\\
	&	=	& \sum_{a,t} (t\tens f_t^{s,a})(s\ra a\tens\delta_{a^{-1}v})\\
	&	=	& \sum_{a,t} t\cdot(s\ra a)\tens\delta_{a^{-1}v}\, f_t^{s,a}((s\ra a)\la (a^{-1}v))\\
	&	=	& \sum_{a,t} t\cdot (s\ra a)\tens\delta_{a^{-1}v}\, f_t^{s,a}((s\la a)^{-1}(s\la v))
\end{eqnarray*}
using the definitions, and ({\ref{7}) in the last step. We multiply both sides from the right by $(e\tens\delta_{u^{-1}v})$ which picks out $a=u$ in the sum. Comparing results,  we see that
\[ \delta_{v,e}e=\sum_t t\cdot(s\ra u) f_t^{s,u}((s\la u)^{-1}(s\la v))\]
By right cancellation in $M$ the basis elements appearing on the right are all distinct. 

We look first at $t=(s\ra u)^{-L}$ which gives $e$. Writing $f=f_t^{s,u}$ for brevity and $L_{s\la u}(f)=f((s\la u)^{-1}(\ ))$, we have
\[ L_{s\la u}(f)=f(s\la(\ ))=\delta_e=\delta_e(s\la(\ ))\]
where the last equality is because $v=e\Leftrightarrow s\la v=e$. To see this, if $v=e$ then $s\la v=e$ by (\ref{8}). Conversely, if $s\la v=e$ then from Lemma~\ref{lemma3.7} we have $v^{-1}=(s\ra v)^{-L}\la(s\la v)^{-1}=(s\ra v)^{-L}\la e=e$. Returning to our displayed equation, we see that $L_{s\la u}(f)=\delta_e$ on $s\la G$. But $s\la G=G$ by Proposition~\ref{prop3.5} hence $f=\delta_{(s\la u)^{-1}}$.

Looking now at $t\ne (s\ra u)^{-L}$, each element of $M$ in our sum occurs just once on the right and not at all on the left. Hence for these $f=f_t^{s,u}$ we have
\[ L_{s\la u}(f)(s\la(\ ))=0\]
and as $s\la G=G$ we conclude that $L_{s\la u}(f)=0$ and hence $f=0$. 

 \endproof

The following applies to the form of $S$ in Lemma~\ref{lemma4.1} which one can also arrive at from other considerations, for example by requiring that $S$ is basis preserving. We continue to state it under the given assumption.

\begin{theorem}
Suppose $M$ in Proposition~3.1 has right-inverses and $G$ is finite. The bicrossproduct $kM\bicross k(G)$ is a Hopf quasigroup iff $M$ is a (two sided inverse property) quasigroup, and  
\begin{equation} \tau(s\ra(t\la u),t\ra u) = (s\la (t\la u))^{-1}\tau(s,t) (s\la (t\la u)) \end{equation}
holds for all $s,t\in M$ and $u\in G$. In this case $S^2=\id$.
\label{theorem}
\end{theorem}
\proof
We suppose that  $S$ has the form found in Lemma~\ref{lemma4.1}. Suppose $kM\bicross k(G)$ is a Hopf quasigroup.  From one of the Hopf quasigroup identities we know that for all $s,t\in M$,
\begin{eqnarray*}
t\tens\delta_e &	=	& S((s\tens\delta_e)\o)((s\tens\delta_e)\t(t\tens\delta_e)) \\
	&	=	& \sum_a S(s\tens\delta_a)((s\ra a\tens\delta_{a^{-1}})(t\tens\delta_e))\\
	&	=	& \sum_a ((s\ra a)^{-L}\tens\delta_{(s\la a)^{-1}})((s\ra a)\cdot t \tens \delta_e) \, \delta_{a^{-1},t\la e}\\
	&	=	& (s^{-L}\tens\delta_e)(s\cdot t\tens\delta_e)\\
	&	=	& s^{-L}\cdot(s\cdot t)\tens\delta_e
\end{eqnarray*}
where only $a$ such that  $a^{-1}=t\la e=e$, i.e. $a=e$ contributes in the sum. So we find that for all $s,t\in M$, $s^{-L}\cdot(s\cdot t)=t$. Similarly on the other side,  we also have
\begin{eqnarray*}
t\tens\delta_e &	=	&((t\tens\delta_e)(s\tens\delta_e)\o)S((s\tens\delta_e)\t)\\
	&	=	& \sum_a (t\cdot s)\cdot s^{-L}\tens \delta_{(s\la a)} \,\delta_{e,s\la a}\delta_{a,s^{-L}\la(s\la a)}\\
	&	=	& (t\cdot s)\cdot s^{-L}\tens \delta_e 
\end{eqnarray*}
where we can replace $s\la a$ by $e$ in view of one of the delta-functions and then only $a=s^{-L}\la e$, i.e. $a=e$ contributes in the sum. So we see that $(t\cdot s)\cdot s^{-L}=t$ for all $s,t\in M$. Hence $M$ is an (inverse property) quasigroup. 

As this is necessary, we will now suppose that $M$ is a quasigroup for the rest of the proof and show that we have a Hopf quasigroup iff the remaining condition displayed in the theorem holds.  Let us see first that it is necessary. We look at one of the Hopf quasigroup identities
\begin{eqnarray*}
t\tens\delta_v	&	=	& (s\tens\delta_e)\o(S((s\tens\delta_e)\t)(t\tens\delta_v))\\
	&	=	& \sum_a (s\tens\delta_a)(S(s\ra a\tens\delta_{a^{-1}})(t\tens\delta_v))\\
	&	=	& \sum_a (s\tens\delta_a)((s^{-L}\tens\delta_{s\la a})(t\tens\delta_v))\\
	&	=	&  s\cdot(s^{-L}\cdot t)\tens\delta_v\,\sum_a \delta_{s\la a,t\la v}\delta_{a,s^{-L}\cdot t \la v}\\
	&=& s\cdot (s^{-L}\cdot t)\tens\delta_v \delta_{s\la((s^{-L}\cdot t)\la v),t\la v}
\end{eqnarray*}
where only $a=s^{-L}.t\la v$ contributes in the sum. We have already dealt with the first tensor factors ($M$ is a quasigroup) and we conclude further that $s\la((s^{-L}\cdot t)\la v)=t\la v$ for all $s,t\in M$ and $v\in G$.  By changing variables we conclude
\begin{equation}\label{eq2} (s\cdot t)\la u=s\la(t\la u)\end{equation}
for all $s,t\in M$ and $u\in G$. Comparing with (\ref{6}) this is equivalent to the condition displayed in the theorem. We have written it entirely in terms of $\tau,\la,\ra$ but clearly it is also equivalent to
\begin{equation}\label{eq3} \tau(s\ra(t\la u),t\ra u) = ((s\cdot t)\la u)^{-1}\tau(s,t) ((s\cdot t)\la u) \end{equation}
 for all $s,t\in M$ and $u\in G$. 
 
 Now suppose that this condition holds and that $M$ is a quasigroup. From a special case of (\ref{eq3}) we see that the condition in Lemma~\ref{lemma3.7} applies for all elements. It remains to verify all of the Hopf quasigroup identities. Thus,
\begin{eqnarray*}
&&\kern-20pt ((t\tens\delta_v)(s\tens\delta_u)\o)S((s\tens\delta_u)\t)\\
&	=	& \sum_a ( t\cdot s\tens \delta_a) \,\delta_{v,s\la a} S(s\ra a\tens\delta_{a^{-1}u})\\
&	=	& \sum_a ( t\cdot s\tens \delta_a) \,\delta_{v,s\la a} ((s\ra u)^{-L}\tens\delta_{((s\ra a)\la(a^{-1}u))^{-1}})\\ 
&	=	& \sum_a ( t\cdot s\tens \delta_a) \,\delta_{v,s\la a} ((s\ra u)^{-L}\tens\delta_{(s\la u)^{-1}(s\la a)})\\
	&	= &\sum_a (t\cdot s)\cdot(s\ra u)^{-L}\tens\delta_{(s\la u)^{-1}v}\delta_{v,s\la a}\delta_{a,(s\ra u)^{-L}\la((s\la u)^{-1}v)}\\
&=& \sum_a (t\cdot s)\cdot(s\ra u)^{-L}\tens\delta_{(s\la u)^{-1}v}\delta_{v,s\la a}\delta_{a,((s\ra u)^{-L}\la(s\la u)^{-1})((s\ra u)^{-L}\ra(s\la u)^{-1})\la v}\\
&=& \sum_a (t\cdot s)(s\ra u)^{-L}\tens\delta_{(s\la u)^{-1}v}\delta_{v,s\la a}\delta_{a, u^{-1}(s^{-L}\la  v)}\\
&=& (t\cdot s)\cdot(s\ra u)^{-L}\tens\delta_{(s\la u)^{-1}v} \delta_{v,s\la(u^{-1}(s^{-L}\la v))}
\end{eqnarray*}
where we used the definitions and (\ref{7}),(\ref{12}) to arrive at the third equality. We then compute the remaining algebra product and replace $s\la a$ by $v$ in some of the expressions in view of $\delta_{v,s\la a}$ to arrive at the 4th equality. Next, we use Lemmas~\ref{lemma3.7} and Lemma~\ref{lemma3.8} and simplify to arrive at the 6th equality.   We now see that only one value of $a$ contributes in the sum to arrive at the final expression.  From the condition on $\tau$ in the theorem, or rather the case of it equivalent to Lemma~\ref{lemma3.7}, we see that
\[ v=s\la(u^{-1}(s^{-L}\la v))\ \Leftrightarrow\ s^{-L}\la v=u^{-1} (s^{-L}\la v)\ \Leftrightarrow u=e.\]
We know here that $s^{-L}\la(\ )$ is bijective by the noted strong form of Proposition~\ref{prop3.5}. Hence we obtain that our original expression
\[ = (t\cdot s)\cdot (s\ra u)^{-L}\tens\delta_{(s\la u)^{-1}v} \delta_{u,e}=(t\cdot s)\cdot s^{-L}\tens\delta_v\delta_{u,e}=t\tens\delta_v\delta_{u,e}\]
as required.

Similarly, doing now the full version,
\begin{eqnarray*}
&&\kern-20pt (s\tens\delta_u)\o(S((s\tens\delta_u)\t)(t\tens\delta_v))\\
	&	=	& \sum_a (s\tens\delta_a)(S(s\ra a\tens\delta_{a^{-1}u})(t\tens\delta_v))\\
	&	=	& \sum_a (s\tens\delta_a)(((s\ra u)^{-L}\tens\delta_{(s\la u)^{-1}(s\la a)})(t\tens\delta_v))\\
	&	=	& \sum_a (s\tens\delta_a)((s\ra u)^{-L}\cdot t\tens\delta_v)\delta_{(s\la u)^{-1}(s\la a),t\la v}\\
	&	=	&  s\cdot((s\ra u)^{-L}\cdot t)\tens\delta_v\,\sum_a \delta_{(s\la u)^{-1}(s\la a),t\la v}\delta_{a,(s\ra u)^{-L}\cdot t \la v}\\
	&=& s\cdot ((s\ra u)^{-L}\cdot t)\tens\delta_v \delta_{s\la(((s\ra u)^{-L}\cdot t)\la v),(s\la u)(t\la v)}
\end{eqnarray*}
by similar computations for the antipode and multiplying out the products. Here only $a=(s\ra u)^{-L}.t\la v$ contributes in the sum and we rearrange the final delta-function accordingly.  But
\[ (s\la((s\ra u)^{-L}\cdot t\la v))=(s\la u)(t\la v)\ \Leftrightarrow\  (s\ra u)^{-L}\cdot t\la v= (s^{-L}\la(s\la u))(((s^{-L}\ra(s\la u))\ra (t\la v))\]
\[ \Leftrightarrow\ (s\ra u)^{-L}\cdot t\la v=u ((s\ra u)^{-L}\cdot t\la v)\ \Leftrightarrow\ u=e\]
using bijectivity of $s^{-L}\la(\ )$, the instance of (\ref{eq2}) in Lemma~\ref{lemma3.7} and (\ref{13}). Hence our original expression
\[ =s\cdot ((s\ra u)^{-L}\cdot t)\tens\delta_v\delta_{u,e}=s\cdot (s^{-L}\cdot t)\tens\delta_v\delta_{u,e}=t\tens\delta_v\delta_{u,e}\]
as required.

We similarly compute
\begin{eqnarray*}
&&\kern-20ptS((s\tens\delta_u)\o)((s\tens\delta_u)\t (t\tens\delta_v))\\	&	=	& \sum_{a}S(s\tens\delta_a)((s\ra a\tens\delta_{a^{-1}u})(t\tens\delta_v)) \\
	&	=	& \sum_{a}((s\ra a)^{-L}\tens\delta_{(s\la a)^{-1}})((s\ra a)\cdot t\tens\delta_v) \,\delta_{a^{-1}u,t\la v} \\
	&	=	& \sum_{a} (s\ra a)^{-L}\cdot((s\ra a)\cdot t)\tens \delta_v \, \delta_{a^{-1}u,t\la v}\delta_{(s\la a)^{-1},((s\ra a)\cdot t)\la v}\\
	&	=	& \sum_{a} t\tens \delta_v \, \delta_{a^{-1}u,t\la v}\delta_{(s\la a)^{-1},((s\ra a)\cdot t)\la v}
\\
&=&t\tens\delta_v \delta_{(t\la v)u,t\la v}=t\tens\delta_v\delta_{u,e}
\end{eqnarray*}
where in one of the delta-functions $(s\ra a)\la a^{-1}=(s\ra a)\cdot t\la v$ is equivalent on applying $(s\ra a)^{-L}\la$ to both sides and using (\ref{eq2}) to $a^{-1}=t\la v$. Hence there is only one value of $a$ in the sum. Finally, 
\begin{eqnarray*}
&&\kern-30pt((t\tens\delta_v)S((s\tens\delta_u)\o)(s\tens\delta_u)\t \\
	&	=	&	\sum_{a}((t\tens\delta_v)S(s\tens\delta_a))(s\ra a\tens\delta_{a^{-1}u})\\
	&	=	&	\sum_{a}((t\tens\delta_v)((s\ra a)^{-L}\tens\delta_{(s\la a)^{-1}}))(s\ra a\tens\delta_{a^{-}u})\\
	&	=	&	\sum_{a} (t\cdot(s\ra a)^{-L}\tens\delta_{(s\la a)^{-1}})(s\ra a\tens\delta_{a^{-1}u}) \, \delta_{v,(s\ra a)^{-L}\la(s\la a)^{-1}} \\
	&	=	& \sum_{a} (t\cdot(s\ra a)^{-L})\cdot(s\ra a)\tens \delta_{a^{-1}u} \, \delta_{v,(s\ra a)^{-L}\la(s\la a)^{-1}}\delta_{(s\la a)^{-1},(s\ra a)\la (a^{-1}u)}\\
	&	=	& \sum_{a} t\tens\delta_{a^{-1}u} \, \delta_{v,a^{-1}}\delta_{(s\la a)^{-1},(s\ra a)\la (a^{-1}u)}\\
	&	=	&  t\tens\delta_{vu} \, \delta_{s\la v^{-1},(s\la u)^{-1}(s\la v^{-1})}=t\tens\delta_{vu} \delta_{s\la u,e}=t\tens\delta_{v} \delta_{u,e}
\end{eqnarray*}
where only $a=v^{-1}$ contributes and we used again that $s\la(\ )$ is bijective.

For the computation of $S^2$ we use Lemma~\ref{lemma3.7}, Lemma~\ref{lemma3.8} and  Proposition~\ref{prop3.6}. \endproof

\begin{remark}\label{remark} We remark that \cite{Beggs03,Zhu} constructed a monoidal category of $M$-graded $G$-modules from the data in Proposition~3.1 and the existence of right inverses. In this case, where there is an obvious multiplicative functor to vectors spaces, one knows by Tannaka-Krein reconstruction   that there
is  a Drinfeld quasi-Hopf algebra\cite{Dri}   generating this category as its modules. This was given in \cite{Zhu} and a close inspection shows that it has a bicrossproduct form. In our dual formulation this category is that of comodules of the bicrossproduct $kM\bicross k(G)$ with $S$ as in Lemma~\ref{lemma4.1} and a Drinfeld Hopf 3-cocycle 
\[   \phi(s\tens\delta_u\tens t\tens\delta_v\tens r\tens\delta_w) = \delta_{u,\tau^{-1}(t,r)}   \delta_{v,e}   \delta_{w,e}\]
\[ \phi^{-1}(s\tens\delta_u\tens t\tens\delta_v\tens r\tens\delta_w)  =\delta_{u,\tau(t, r)}   \delta_{v,e}   \delta_{w,e}\]
making it into a coquasi-Hopf algebra.  Here\cite{Ma:book} the algebra product is associative up to conjugation by $\phi$ in a convolution sense 
\[\sum \phi(h\o\tens g\o\tens f\o)(h\t g\t) f\t=\sum h\o (g\o f\o)\phi(h\t\tens g\t\tens f\t)\]
for all $h,g,f$ and $\phi$ is invertible in the same convolution sense. The coquasi-Hopf structure in our case
is easily be verified by direct computation and does not
require the two further conditions in our theorem, i.e. the bicrossproduct being a Hopf quasigroup is 
strictly stronger. 

Incidentally, this remark means that there is a coquasi-Hopf algebra and monoidal category associated to any any finite left quasigroup $M$ with right inverses as every such can be expressed as a coset construction\cite{Smi:book}. This is such that, in view of (\ref{5}), the algebra $kM$ is associative in the category. It is not clear if our additional requirements for a Hopf quasigroup can ever be satisfied, starting now from $M$ an IP quasigroup. One can show that the following special case cannot.
 \end{remark}

\begin{corollary}\label{corol4.4} We obtain a bicrossproduct Hopf quasigroup in the context of Theorem~\ref{theorem} if $M$ is an (inverse property) quasigroup and
\[ \tau(s, t\ra u)=\tau(s,t),\quad \tau(s\ra(s^{-L}\la u),t)=u^{-1}\tau(s,t)u\]
for all $s,t\in M$ and $u\in G$.
\end{corollary}
\proof If we suppose that the first of these conditions then the requirement on $\tau$ in Theorem~\ref{theorem} becomes on a change of variables $t\la u$ to $u$ (possible by $M$ a quasigroup and the strong form of Proposition~\ref{prop3.5}), 
\[ \tau(s\ra u,t)=(s\la u)^{-1}\tau(s,t)(s\la u).\]
This is equivalent to the second condition stated, on a further change of variables $s\la u$ to $u$.  \endproof

The first condition in the corollary says that the second argument of $\tau$ is constant on orbits of $\ra$. The 2nd condition says that the function in its first argument 
essentially intertwines  $\ra$ with the adjoint action (cf. a crossed module $\tau:M\to G$) except that it is  twisted by $\la$.  This suggests to further simplify our search by focusing on the special case where $\la$ is trivial. We denote by $[t]$ the orbit label or equivalence class of $t$ under the remaining action $\ra$. Note that if $\la$ is trivial then the condition on $\tau$ in Theorem~\ref{theorem} is already included as (\ref{6}) in Proposition~3.1, i.e. there is no additional constraint in this case other than $M$ an IP quasigroup. However, we still have to solve for this data and we will do so in the special case of the corollary.

\begin{corollary}\label{corol4.5} Suppose that $\la$ is trivial. Then the data for constructing a bicrossproduct Hopf  
quasigroup in the special case of Corollary~\ref{corol4.4} become $G$ finite and $\ra,\tau$ such that
\begin{enumerate}
\item[(i)] $\tau(s,t)=\tau(s,[t])$  (depends only on $[t]$)
\item[(ii)] $u^{-1} \tau(s,t)u=\tau(s\ra u, \ra u)$   (covariance condition)
\item[(iii)] $\tau(s,[t])\tau(s\cdot t,[r])=\tau(s,[t\cdot r])\tau(t,[r])$ \label{2-cocycle}(2-cocycle  
condition).
\end{enumerate}
Here we require that $G$ acts on $M$ by an action $\ra$ respecting its structure as in (\ref{3})-(\ref{4}), (\ref{8}) and that $M$ is an IP quasigroup and quasiassociative in the sense (\ref{5}). 
\end{corollary}
\proof The condition (\ref{eq1}) coincides with (\ref{6}) and with (ii) when $\la$ is trivial, (\ref{7}) is empty,  while (\ref{2}) simplifies to (iii) on use of (i). We still require (\ref{3})-(\ref{4}), (\ref{5}) and (\ref{8}) as stated. Clearly $s\la(\ )$ is bijective as it is the identity so we have right inverses by Proposition~\ref{prop3.5}. \endproof

In this case $X$ has a semidirect product form with cross relations $su=u (s\ra u)$ and relates to the product of $M$ by $st=\tau(s,t)s\cdot t$ as before. Its structure can be recovered given the stated data from
\[  (u,s)(v,t)=(u\tau(s,[t])v,[(s\ra v)\cdot t])\]
as a special case of Proposition~\ref{Beggs}. Also, if $G$ is Abelian then the first two conditions say that $\tau(s,t)=\tau([s],[t])$ depends for both of its arguments only on the orbits in $M$ under $G$. 

For an example we let $M=G_\O$ the octonion quasigroup which we will take in the binary-vector form of the octonions introduced in \cite{AlbMa}. The group here consists of elements $\{\pm e_{\vec a}\}$, where $\vec a\in \Z_2^3$ is a 3-tuple with values in $\{0,1\}$, sitting inside the octonion algebra with product $e_{\vec a}\cdot e_{\vec b}=F(\vec a,\vec b)e_{\vec a+\vec b}$ in terms of component-wise addition. The signs here are given by\cite{AlbMa}
\[ F(\vec a,\vec b)=\left(
\begin{array}{cccccccc}
 1 & 1 & 1 & 1 & 1 & 1 & 1 & 1 \\
 1 & -1 & 1 & -1 & 1 & -1 & 1 & -1 \\
 1 & -1 & -1 & 1 & 1 & -1 & -1 & 1 \\
 1 & 1 & -1 & -1 & -1 & -1 & 1 & 1 \\
 1 & -1 & -1 & 1 & -1 & 1 & 1 & -1 \\
 1 & 1 & 1 & 1 & -1 & -1 & -1 & -1 \\
 1 & -1 & 1 & -1 & -1 & 1 & -1 & 1 \\
 1 & 1 & -1 & -1 & 1 & 1 & -1 & -1
\end{array}
\right)\]
in binary basis order $000,001,010,011,\cdots,111$.  The quasigroup in this form, like the octonion algebra, is quasiassociative in the sense\cite{AlbMa}
\[ (e_{\vec a}\cdot e_{\vec b})\cdot e_{\vec c}=\phi(\vec a,\vec b,\vec c)e_{\vec a}\cdot (e_{\vec b}\cdot e_{\vec c}),\quad \phi(\vec a,\vec b,\vec c)=(-1)^{|\vec a,\vec b,\vec c|}\]
extended to signs, where we use the determinant of the matrix formed by the three vectors (in other words $\phi$ is -1 if and only if the three vectors are linearly independent in $\Z_2^3$ as a vector space over $\Z_2$).

\begin{example} Let $X=\Z_2^3\ltimes Cl_3$ be the order 128 non-Abelian group with generators $\pm e_i,g_i$, $i=1,2,3$ and relations
\[ e_ie_j=\begin{cases}-e_j e_i &{\rm if\ }i\ne i\cr -1& {\rm if\ }i=j\end{cases},\quad e_ig_j=\begin{cases}g_j e_i &{\rm if\ }i\ne i\cr -g_je_i& {\rm if\ }i=j\end{cases},\quad g_ig_j=g_jg_i,\quad g_i^2=1\]
and let $G=\Z_2^3$ as generated by $\{g_i\}$. Then the transversal $M\subset X$ labeled by $\pm e_{\vec a}$ and consisting of
\[ \begin{array}{clcl}
e_{000} & = 1,&
e_{001}&=e_3\\
e_{010}&=e_2,&
e_{011}&=-g_1 e_2 e_3\\
e_{100}&=e_1,&
e_{101}&=-g_2e_1e_3\\
e_{110}&=-g_3 e_1e_2,\quad&
e_{111}&=g_1g_2g_3 e_1e_2e_3\end{array}\]
extended to signs, acquires the structure of the octonion quasigroup $G_\O$. Moreover, the conditions of Theorem~\ref{theorem} hold and we have a Hopf quasigroup $k G_\O\lcocross k(\Z_2^3)$. 
\end{example}
\proof This is constructed using Corollary~\ref{corol4.5} above. We know as in \cite{Klim09} that $G_\O$ has an action $\ra$ of the group $\Z_2^3$ given by $e_{\vec a}\ra g^{\vec b}=e_{\vec a}(-1)^{\vec a\cdot \vec b}$ where $g^{\vec b}=g_1^{b_1}g_2^{b_2}g_3^{b_3}$ is the group written multiplicatively. We used the vector space dot product over $\Z_2$. Explicitly, $e_{\vec a}\ra g_i=(-1)^{a_i}e_{\vec a}$. We know that this action respects the product. Next, we define
\[ \tau(\pm e_{\vec a},\pm e_{\vec b})=g^{\vec a\times\vec b}\equiv g_1^{a_2b_3-a_3b_2}g_2^{a_3b_1-a_1b_3}g_3^{a_1b_2-a_2b_1}\]
where the vector space cross product in $\Z_2^3$ is viewed multiplicatively via $g$ as shown explicitly. This obeys the 2-cocycle condition (\ref{2-cocycle}) in Corollary~\ref{corol4.5} due to linearity over $\Z_2$ of the cross product.  We also require that the quasigroup obeys
\[ (e_{\vec a}\cdot e_{\vec b})\cdot e_{\vec c}=e_a\ra\tau(e_{\vec b},e_{\vec c})\, \cdot(e_{\vec b}\cdot e_{\vec c})\]
which obtains because $\vec a\cdot(\vec b\times\vec c)=|\vec a,\vec b,\vec c|$. Similarly when there are $\pm$ signs. Finally, our special conditions for Theorem~\ref{theorem} hold because $\tau$ does not depend on the sign of its arguments and the orbit of $e_{\vec a}$ under $\ra$ is $\{\pm e_{\vec a}\}$ (the group $G$ being abelian, we require that both arguments depend only through the orbits). Hence we have all the data for a quasigroup double cross product to yield a group $X=G\ltimes G_\O$. It remains to determine what this group is. It contains $G=\Z_2^3$ as a subgroup, and cross relations $e_i g_j=g_j e_i\ra g_j $ where $e_1=e_{100}$, $e_2=e_{010}$ and $e_3=e_{001}$ are (at this stage) elements of the quasigroup $G_\O$. This gives the cross-relations stated. We also have $e_{\vec a} e_{\vec b}=\tau(e_{\vec a},e_{\vec b})e_{\vec a}\cdot e_{\vec b}=g^{\vec a\times\vec b} F(\vec a,\vec b)e_{\vec a+\vec b}$ for the product in $X$ in terms of that in $G_\O$. Thus
\[ e_1e_2=g_3F(100,010)e_{110}=-g_3 e_{110}=-e_2e_1\]
\[e_1e_3=g_2F(100,001)e_{101}=-g_2e_{101}=-e_3 e_1\]
\[ e_2e_3=g_1F(010,001)e_{011}=-g_1e_{011}=-e_3e_1,\quad e_i^2=-1\]
where $F(\vec a,\vec b)=-F(\vec b,\vec a)$ when $\vec a,\vec b,\vec a+\vec b\ne0$ (the altercommutativty of octonions, see\cite{AlbMa}). This gives the relations of $X$ in terms of the $\{e_i\}$ regarded now as generators of $X$. These relations and those of $\Z_2^3$ provide for a basis $\{\pm e_1^{a_1}e_2^{a_2}e_3^{a_3} g_1^{b_1}g_2^{b_2}g_3^{b_3}\}$ which has order 128, hence these are all the relations. We also see in these calculations how the products are related to elements of $G_\O$ and rearrange them to obtain the image of most of them in $X$. We similarly compute 
\[ e_1e_2e_3=-g_3e_{110}e_{001}=-g_3g_1g_2F(110,001)e_{111}=g_1g_2g_3e_{111}\]
to obtain the last element $e_{111}$ of the transversal. \endproof

We find that the group $X$ here is a semidirect product by $\Z_2^3$ of the `Clifford group' $Cl_3$ generated by the $\pm e_i$. This is the set of signed monomials of these generators in the Clifford algebra in three dimensions (generalizing the way in which the quaternion group is defined from the quaternion algebra) and easily seen to form a group. The additional information provided by the transversal provides the quasigroup structure on the left coset space $G\backslash X$ according to our results above. For example,
\[ e_{110}e_{001}=-g_3e_1e_2e_3=-g_1g_2e_{111}\] 
induces $Ge_{110}\cdot G e_{001}=G (-e_{111})$ at the level of cosets. In this way one can verify all the signs in the table of $F$ for the $G_\O$ product as a useful check of all of our theory. The additional signs beyond those from the group $Cl_3$ come from moving all the $\{g_i\}$ to the far left where it is absorbed by $G$ in the coset.  We also obtain, of course, a new Hopf quasigroup. Its dual is a Hopf coquasigroup $k(G_\O)\rtimes k\Z_2^3$ with  structure is similar to that of $k[S^7]\rtimes k\Z_2^3$ in \cite{Klim09} and could be obtained in a similar way as there. However, we have provided now a bicrossproduct point of view on it. By Remark~4.3 we also have the bicrossproduct as a coquasi-Hopf algebra and a monoidal category associated to $G_\O$ in this way 

It remains to find interesting examples of the most general $k M\bicross k(G)$ form with all of $\la,\ra,\tau$ nontrivial.

\section*{Acknowledgements} We would like to thank the organisers for a most stimulating conference, and in particular J.D.H. Smith for some discussion.

\end{document}